\documentclass{article}
 \usepackage{amsmath}
  \usepackage{amsthm}
 \usepackage{amssymb}

      \usepackage[english]{babel}
         \usepackage{wrapfig}
           \usepackage{graphicx}
             \usepackage{geometry}
               \usepackage[rightcaption]{sidecap}
                  \usepackage{tikz}
                   \usepackage{bbm}
						\usepackage{url}
\usepackage{amsfonts}
\usepackage{marvosym}
\usepackage{wasysym}
\usepackage{pifont}
\usepackage{mathrsfs}
\DeclareMathAlphabet{\mathpzc}{OT1}{pzc}{m}{it}

\allowdisplaybreaks[2]

\title{An Obstacle for Higher Regularity\\ of Geodesics in the Space of K\"ahler Potentials}
\author{Jingchen Hu}

\newcommand{\Anforderung}{{6}}
\newcommand{\Add}{4}

\newcommand{\skype}{\mathcal{H}}
\newcommand{\Cglatt}{C^{\infty}}
\newcommand{\sigmaone}{\sigma_1}
\newcommand{\sigmatwo}{\sigma_2}
\newcommand{\sS}{\mathcal{S}}
\newcommand{\Real}{{\text{Re\ }}}
\newtheorem{problem}{Problem}[section]
\newtheorem{question}{Question}[section]

\newcommand{\dD}{D}

\newtheorem{notation}{Notation}

\newtheorem{Beispiel}{Example}[section]

\newtheorem{Satz}{Theorem}[section]

\newtheorem{remark}{Remark}[section]

\newcommand{\Gleichung}{equation}

\newcommand{\vv}{{\bf v}}
\newcommand{\vu}{{\bf u}}

\numberwithin{equation}{section}

\newcommand{\ER}{\mathbb{R}}
\newcommand{\EZ}{\mathbb{Z}}
\newcommand{\EC}{\mathbb{C}}

\newcommand{\Torus}{\mathcal{T}}

\newcommand{\ddbar}{\partial\overline\partial}

\newcommand{\spacecdot}{{\ \cdot\ }}

\newcommand{\EchtZentralFaser}{{[0,1]\times \{0\}}}
\newcommand{\ZentralFaser}{\EchtZentralFaser}
\newcommand{\LV}{{\triangle_T}}
\newcommand{\NV}{{\nabla_T}}
\newcommand{\xyVierKlein}{{O(x^4+y^4)}}

\newcommand{\HappytoN}{{\mathcal{P}^{2n}}}

\newcommand{\aoneone}{{a}}
\newcommand{\atwotwo}{{b}}

\newcommand{\HappytoNE}{\mathcal{P}^{2n}_E}

\newcommand{\Lappen}{\LV}
\newcommand{\Naben}{\NV}%
\newcommand{\MEA}{{\mathcal{E}_A}}
\newcommand{\MSA}{{\mathcal{S}_A}}
\newcommand{\Endo}{\text{End}}
\newcommand{\Auto}{\text{Aut}}
\newcommand{\MU}{{\mathcal{U}}}
\newcommand{\Boost}{{\mathcal{B}}}
\newcommand{\MtD}{\widetilde{\mathcal{D}}}
\newcommand{\MD}{{\mathcal{D}}}

\newcommand{\LP}{\mathcal{U}_{LP}}
\newcommand{\ReducetoBase}{\mathcal{R}}

\begin{document}
\maketitle

\begin{abstract} In this paper we address the following question regarding the regularity of geodesics in the space of K\"ahler potentials. Given a geodesic which is highly regular, and has smooth boundary value, can we expect that it is actually smooth? We construct an example to show that the answer to the question is ``no".
\end{abstract}

\section{Introduction}\label{IntroductionSection}
  Given $(V,\omega_0)$, a smooth  K\"ahler manifold, we consider the space of K\"ahler potentials 
\[\skype=\{\phi\in\Cglatt(V)|\omega_0+\sqrt{-1}\ddbar\phi>0\}.\]
At any point $\varphi\in\skype$, the tangent space $T_\varphi\skype$ can be identified with $\Cglatt(V)$, and following Mabuchi \cite{Mabuchi}, we define the following Riemannian metric in $\skype$, for $\psi_1,\psi_2\in T_\varphi\skype$,
\[<\psi_1,\psi_2>_\varphi=\int_V\psi_1\psi_2(\omega_0+\sqrt{-1}\ddbar\varphi)^n.\]
With this metric, the energy of a differentiable curve $\varphi:[0,1]\rightarrow \skype$ is 
\[\int_{0}^1\int_V\left(\frac{d\varphi}{dt}\right)^2(\omega_0+\sqrt{-1}\ddbar\varphi(t))^n dt.\]
Then the geodesic equation for a smooth curve is 
\begin{equation}\varphi_{tt}-g_{\varphi}^{i\overline j}\varphi_{ti}\varphi_{t\overline j}=0.\label{SmoothGeodesicEquation}
\end{equation}
As discovered by Semmes\cite{Semmes} and Donaldson\cite{DonaldsonHolomorphicDiscs}, (\ref{SmoothGeodesicEquation}) can be written as a homogenous complex Monge-Amp\`ere equation. Denote 
\[\sS=\{\tau=t+\sqrt{-1}\theta\in \EC| \ 0\leq t\leq 1\}.\] Then we can consider a curve $\varphi\in C^1([0,1]; \skype)$ as a function defined on  $\sS\times V$, by letting 
\[\Phi(\tau,\ast)=\varphi(\Real\tau).\]
Let $\pi$ be the projection  $\sS\times V\rightarrow V$, and $\Omega_0=\pi^\ast(\omega_0)$. Then for a smooth curve $\varphi:[0,1]\rightarrow \skype$ satisfying  (\ref{SmoothGeodesicEquation})  is equivalent to the corresponding $\Phi$ satisfying  
\[(\Omega_0+\sqrt{-1}\ddbar\Phi)^{n+1}=0.\]

So the problem of finding a geodesic, in $\skype$, connecting $0$ and $\varphi$, can be related to solving the following Dirichlet  problem for the homogenous complex Monge-Amp\`ere equation on $\sS\times V$.\footnote{In this paper, we abbreviate ``homogenous complex Monge-Amp\`ere equation" as ``HCMA equation''.}

\begin{problem}{[Dirichlet Problem for the HCMA Equation on $\sS\times V$]}\label{Problem:Dirichlet_HCMA_Strip}\\
Given $\varphi\in\skype$,  find a $\Phi\in C^2(\sS\times V)$, so that
\begin{align}
(\Omega_0+\sqrt{-1}\ddbar\Phi)^{n+1}=0,\ \ &\text{ in }\sS\times V;\\
\Phi(\tau, \ast)=0, \ \ &\text{ for }\Real\tau=0;\\
\Phi(\tau, \ast)=\varphi, \ \ &\text{ for }\Real\tau=1;\\
\partial_\theta\Phi=0, \ \ &\text{ in } \sS\times V;\\
\omega_0+\sqrt{-1}\ddbar\Phi(\tau,\ast)\geq 0, \ \ &\text{ for }\tau\in\sS.
\end{align}
\end{problem}
Due to the work of
Chen\cite{ChenPrinceton}, 
Chu-Tossati-Weinkove\cite{Jianchun}, 
 the problem above always has a $C^{1,1}$ solution. But in general a solution may not correspond to a curve in $\skype$. First $\Phi(\tau,\ast)$ may not be $C^{\infty}$. As shown by Liz Vivas-Lempert-Darvas\cite{LempertLizVivas}\cite{LempertDarvas}\cite{DarvasMorse}, $C^{1,1}$ is the optimal global regularity for general $\varphi$. Second, it is expected that $\omega_0+\sqrt{-1}\ddbar\Phi(\tau,\ast)$ may degenerate for some $\tau\in \sS$, which makes $\Phi(\tau,\ast)\notin\skype$. However, we can consider a solution to Problem \ref{Problem:Dirichlet_HCMA_Strip} as a weak or generalized geodesic.

 In this paper, if $\Phi$,  a solution to Problem \ref{Problem:Dirichlet_HCMA_Strip}, satisfies
\[\omega_0+\sqrt{-1}\ddbar\Phi(\tau,\ast)> 0, \ \ \ \ \ \ \text{ for }\tau\in\sS,\]
we say $\Phi$ is a {\bf non-degenerat}e geodesic. And if $\Phi\in C^k(\sS\times V)$ we say that $\Phi$ is a {\bf $C^k$ geodesic}. Similarly, if $\Phi\in \Cglatt(\sS\times V)$, we say $\Phi$ is a {\bf smooth geodesic}.

A problem analogous to Problem \ref{Problem:Dirichlet_HCMA_Strip} is the following Dirichlet problem for the HCMA equation on the product of a disc and a manifold.

\begin{problem}{[Dirichlet  Problem for the HCMA Equation on $D\times V$]}\label{Problem:Dirichlet_HCMA_Disc}\\Let $D$ be the unit disc in the complex plane.
Given $(V,\omega_0)$, a smooth K\"ahler manifold,  and $F\in C^{\infty}(\partial D\times V)$, satisfying 
\[\omega_0+\sqrt{-1}\ddbar F(\tau,\ast)>0,\ \ \ \ \text{ for }\tau\in\partial\dD,\]
 find $\Phi\in C^{2}(D\times V)$, satisfying
\begin{align}
(\Omega_0+\sqrt{-1}\ddbar\Phi)^{n+1}=0,\ \ &\text{ in }\dD\times V;\\
\Phi(\tau, \ast)=F, \ \ &\text{ for }\tau\in\partial \dD;\\
\omega_0+\sqrt{-1}\ddbar\Phi(\tau,\ast)> 0, \ \ &\text{ for }\tau\in\dD.
\end{align}
\end{problem}
In \cite{DonaldsonHolomorphicDiscs}, by relating the  Dirichlet problem for the HCMA equation to the existence and stability of a family of holomorphic discs with boundaries attached to a totally real submanifold, Donaldson proved
\begin{Satz}{(Donaldson, \cite{DonaldsonHolomorphicDiscs})}\label{Satz:DonaldsonDiscStableSolvability}
The set of boundary values $F$, for which Problem \ref{Problem:Dirichlet_HCMA_Disc} has a smooth solution, is an open set in $\Cglatt(\partial\dD\times V)$ with respect to the $C^3$ topology.
\end{Satz}
Another consequence of applying the technique of \cite{DonaldsonHolomorphicDiscs} is:
\begin{Satz}\label{Satz:DonaldsonBootStrap}
Given $F\in \Cglatt(\partial \dD\times V)$, if  $\Phi\in C^3(D\times V)$ is a solution to Problem \ref{Problem:Dirichlet_HCMA_Disc}, then $\Phi\in\Cglatt(\dD\times V)$.
\end{Satz}

In \cite{CFH}, by partially generalizing the technique of \cite{DonaldsonHolomorphicDiscs} to the strip case, we proved
\begin{Satz}{(Chen-Feldman-Hu,\ Theorem 1.2 \cite{CFH})}\label{Satz:CFH}\\
For any fixed $k>4$, there exists a $\delta>0$, so that if $\varphi\in \Cglatt(V)$ satisfies $|\varphi|_k<\delta$, the geodesic connecting $0$ and $\varphi$ is $C^4$ and non-degenerate.
\end{Satz}
\begin{notation}In \cite{CFH} and this paper, we use $|\cdot|_k$ to denote the $C^k$ norm.
\end{notation}
 Comparing Theorem \ref{Satz:DonaldsonDiscStableSolvability} and  Theorem \ref{Satz:DonaldsonBootStrap} with our Theorem \ref{Satz:CFH}, it is natural to ask  if any  analogues of Theorem \ref{Satz:DonaldsonDiscStableSolvability} and Theorem \ref{Satz:DonaldsonBootStrap} are available when the disc is replaced by a strip? For example we ask:
\begin{question}
\label{Vermutung:KurzGlatt}
Given $(V, \omega_0)$ a K\"ahler manifold, can we find a $B=B(V, \omega_0)$ and $\varepsilon=\varepsilon(V,\omega_0), $ so that if  $\varphi\in C^{\infty}(V)$ satisfies
\[|\varphi|_B\leq \varepsilon,\]
then there exists a  non-degenerate smooth  geodesic $\Phi$ connecting $0$ and $\varphi$?
\end{question}

{\bf\flushleft Main Result:}
We construct the following example to show that the answer to Question \ref{Vermutung:KurzGlatt} is ``no".

\begin{Beispiel}\label{Beispiel:Kurz}
On the torus $\Torus=S^1\times S^1$, with flat background metric 
\[\omega_0=\frac{1}{2}dx\wedge dy,\]
there exists a sequence of analytic functions $\varphi_k$, for k=1,2,... ,  so that
\[|\varphi_k|_B\rightarrow 0, \text{ as }k\rightarrow \infty,\]
for any fixed $B$, but none of $\varphi_k$ can be connected with $0$ by a smooth non-degenerate geodesic.
\end{Beispiel}

\begin{remark}
In Example \ref{Beispiel:Kurz}, let us denote the geodesic connecting $0$ and $\varphi_k$ by $\Phi_k$, then by Theorem 1.8 of \cite{CFH}, we know for $k$ big enough, $\Phi_k$ is  non-degenerate and $\Phi_k$ gets more and more regular as $k\rightarrow \infty$. But the example shows, none of $\Phi_k$ can be $\Cglatt$. 
\end{remark}

In Section \ref{ODEsystemsforderivatives20190105}, we study the ODE systems satisfied by the derivatives of a geodesic. In Subsection \ref{ODESystemofSecondOrderDerivativesofGeodesic},  the ODE system for the second derivatives is studied and some properties, which will be used, are listed, however proofs of some properties are postponed to the Appendix \ref{Solvability_of_2nd_Jets}. In Subsection \ref{ODESystemofHigherOrderDerivativesofGeodesic}, assuming a geodesic $\Phi$ is regular enough, we analyze the ODE system for the higher derivatives, and get Theorem \ref{Satz:Paaren}, which says that, in certain situations, the derivatives of $\Phi(0,*)$ and $\Phi(1,*)$ should satisfy a condition. Based on this, we construct the claimed example in Section \ref{KurzUnGlatt}.

\section{ODE Systems for the Derivatives of Geodesics}\label{ODEsystemsforderivatives20190105}

Our background manifold will be the torus
\[\Torus=[-\pi,\pi]\times[-\pi,\pi]\slash\sim,\]
where
\[(-\pi,y)\sim(\pi,y), \text{ for any } y\in[-\pi,\pi],\]
and
\[(x, -\pi)\sim(x,\pi), \text{ for any } x\in[-\pi,\pi].\]
The coordinates on $\Torus$ will be denoted by $(x,y)$. The background K\"ahler form $\omega_0$ is \[\frac{dx\wedge dy}{2}.\]

In this and the next section, we consider a curve in $\skype$ as a function defined on $[0,1]\times \Torus$. Then, with our background K\"ahler form, 
the geodesic equation can be simplified as
\begin{\Gleichung}
(1+\LV\Phi)\Phi_{tt}=|\NV\Phi_t|^2.   				\label{Gleichung:GeodesicEquationonToruswithFlatBackGround}
\end{\Gleichung}
\begin{notation}
Here, $\LV=\partial_x^2+\partial_y^2$ and $\NV\spacecdot=(\partial_x\spacecdot,\ \partial_y\spacecdot)$.
\end{notation}

Given two potentials $\varphi_0$, $\varphi_1$, which are even both in $x$ and $y$, i.e. satisfying
\[\varphi_0(x, y)=\varphi_0(-x, y), \ \varphi_0(x, y)=\varphi_0(x, -y), \]
\[\varphi_1(x, y)=\varphi_1(-x, y), \ \varphi_1(x, y)=\varphi_1(x, -y), \]
 according to the uniqueness theorem for the boundary value problem: Corollary 7 of \cite{DonaldsonSymmetricSpace},  if there is a geodesic $\Phi$ connecting $\varphi_0$ and $\varphi_1$, $\Phi$ should also be even in both $x$ and $y$, i.e. satisfies
\[\Phi(t,x, y)=\Phi(t,-x, y), \ \Phi(t,x, y)=\Phi(t,x, -y), \]
for any $t\in[0,1].$
For such $\Phi$, all odd order $x$ or $y$-derivatives vanish along $\ZentralFaser$, i.e.
\begin{equation}D_t^k D_x^i D_y^j\Phi=0, \text{ on } \ZentralFaser, \text{ if }2\nmid i \text{, or }2\nmid j.\label{SymmetryJan052019}
\end{equation}

For the convenience of the presentation, we introduce the following notations:
\begin{notation}
 For a function $f$ on $\Torus$, we call
\[\sum_{0\leq i, j }^{i+j\leq k}(D_x^i D_y^j f)(p)x^iy^j\] the $k$-jet of $f$ at  the point $p$.
For a function $F$ defined on $[0,1]\times \Torus$, we call
\[\sum_{0\leq i, j }^{i+j\leq k}(D_x^i D_y^j f)(t, p)x^iy^j\]the $\Torus$-directional $k$-jet of $F$ at $(t, p)$, for $0\leq t\leq 1$.

\end{notation}


\subsection{Second Order Derivatives}\label{ODESystemofSecondOrderDerivativesofGeodesic}
Suppose we have a non-degenerate geodesic $\Phi\in C^{\Anforderung}([0,1]\times V)$, 
connecting $\varphi_0$ and $\varphi_1$, which are even in both $x$ and $y$. Then $\Phi$ should  also be even in both $x$ and $y$, so
\[\Phi_{tt}=\frac{|\NV\Phi_t|^2}{1+\LV\Phi}=0,\text{ on } \ZentralFaser.\]
In addition, we assume $\varphi_0(0)=\varphi_1(0)=0$, hence
\[\Phi\equiv0, \ \text{ on }\ZentralFaser.\]

\begin{notation}
Prime `` $'$ " will be used to denote $\frac{d}{dt}$.
\end{notation}

Along $\ZentralFaser$, let the second order Taylor expansion of $\Phi$, with respect to $x$ and $y$,  be
\begin{equation}\Phi=\aoneone x^2+\atwotwo y^2+\xyVierKlein.  \label{2JetsJan052019}
\end{equation}
Above, $\aoneone $ and $\atwotwo $ are $C^2$ functions of $t$.
By simply plugging (\ref{2JetsJan052019}) into equation (\ref{Gleichung:GeodesicEquationonToruswithFlatBackGround}), we get, on $[0,1]$,
\begin{align}\renewcommand\arraystretch{1.3}
\aoneone ''=\frac{4(\aoneone ')^2}{1+2\aoneone +2\atwotwo }	;		\label{Evolutionofa11}\\
\atwotwo ''=\frac{4(\atwotwo ')^2}{1+2\aoneone +2\atwotwo }	.		\label{Evolutionofa22}
\end{align}
We denote 
\begin{align}\sigmatwo=\frac{\aoneone ' \atwotwo '}{(1+2\aoneone +2\atwotwo )^2},\ \ \ \ \ 
\sigmaone=\frac{\aoneone '+\atwotwo '}{1+2\aoneone +2\atwotwo }.		\label{SigmaonetwoJan032019}
\end{align}
Taking derivative of (\ref{SigmaonetwoJan032019}) and using (\ref{Evolutionofa11}) (\ref{Evolutionofa22}), we get
\begin{\Gleichung}\sigmatwo'=0.\label{InvarianceofBifurcatingStationaryLocus}
\end{\Gleichung}
\begin{\Gleichung}\sigmaone'=\frac{2(\aoneone '-\atwotwo ')^2}{(1+2\aoneone +2\atwotwo )^2}=2\sigmaone^2-8\sigmatwo\geq 0.		\label{EntwicklungPfad}
\end{\Gleichung}

So $\sigmatwo$ is a constant on $[0,1]$. Our examples will be constructed in the situation of $\sigmatwo<0$. 
We assert that  $\sigmatwo<0$, if and only if 
\[(\aoneone (1)-\aoneone (0))\cdot(\atwotwo (1)-\atwotwo (0))<0.\]
 This is proved by Theorem \ref{Satz:SecondJetsLemmaA} of Appendix \ref{Solvability_of_2nd_Jets}.

In the following, we assume that $\sigma_2<0$ and denote $\sigmatwo\equiv-\epsilon^2$, for some constant $\epsilon>0$. We can also assume that $a'>0$
and let
\begin{\Gleichung}\frac{\aoneone '}{1+2\aoneone +2\atwotwo }=\epsilon A, \ \ \ \ \ \frac{\atwotwo '}{1+2\aoneone +2\atwotwo }=\frac{-\epsilon}{A},\label{Gleichung:Derivativeofa11a22overVolumeform}
\end{\Gleichung}
with $A$ a positive real valued function of $t$, if $a'<0$ we can interchange $x$ and $y$. 
Plugging (\ref{Gleichung:Derivativeofa11a22overVolumeform}) into the definition of $\sigma_1$ in (\ref{SigmaonetwoJan032019}), taking $\frac{d}{dt}$-derivative of $\sigma_1$, then using  (\ref{EntwicklungPfad}), we get
\begin{equation}\epsilon(1+\frac{1}{A^2})A'=(\epsilon A-\frac{\epsilon}{A})'=\sigmaone'=2\left(\frac{\aoneone '-\atwotwo '}{1+2\aoneone +2\atwotwo }\right)^2=2\epsilon^2(A+\frac{1}{A})^2,\label{Plugginginsigma1_20180803}
\end{equation}
which can be reduced to 
\begin{equation}A'=2\epsilon (A^2+1).\label{ODE_of_A_hyperbolic}
\end{equation}

Then taking $\frac{d}{dt}$-derivative of (\ref{Gleichung:Derivativeofa11a22overVolumeform}) gives
\begin{equation}\label{Gleichung:secondDerivativeofa11a22overVolumeform}
\frac{\aoneone ''}{1+2\aoneone +2\atwotwo }=4\epsilon^2A^2,\ \ \ \ \ \frac{\atwotwo ''}{1+2\aoneone +2\atwotwo }=4\epsilon^2\frac{1}{A^2}.
\end{equation}

In Appendix \ref{Solvability_of_2nd_Jets}, Theorem \ref{Satz:SecondJetsLemmaA} shows that,  if  $\aoneone , \ \atwotwo \in C^2([0,1])$ satisfy equation (\ref{Evolutionofa11}) (\ref{Evolutionofa22}), then they are uniquely determined by their boundary values, i.e. determined by $\aoneone (0),$ $ \aoneone (1),$ $ \atwotwo (0)$, $\atwotwo (1)$. This means that if there is a $C^\Anforderung$ non-degenerate geodesic connecting $\varphi_0$ and $\varphi_1$ then, along $\ZentralFaser$, all the second derivatives of $\Phi$ are determined by the 2-jets of $\varphi_0$ and $\varphi_1$ at $0\in\Torus$.
And the $\epsilon$ above is determined by the boundary value of $\aoneone, \atwotwo$, by formula (\ref{LemmaA_Formula_for_epsilon}).


\subsection{Higher Derivatives} \label{ODESystemofHigherOrderDerivativesofGeodesic}
Having studied the second derivatives of  $\Phi$, we turn to discuss the higher derivatives of $\Phi$ along  $\ZentralFaser$. We still assume the geodesic $\Phi$ is even in both $x$ and $y$. And to guarantee the validity of the following computation we assume that, for some $n\geq 2$, $\Phi\in C^{2n+4}([0,1]\times \Torus)$. 

Along $\ZentralFaser$, let the Taylor expansion of $\Phi$, with respect to $x,y$,   be
\begin{equation}
\Phi=L(t, x,y )+P(t,x,y)+R(t,x,y),\label{RevisedTaylorDec302018}
\end{equation}
where, for each fixed $t$,  $L$ is a polynomial in $x,y$ of degree smaller than $2n$,  $P$ is homogenous of  degree $2n$, and for some constant $C$, $$R(t,x,y)\leq C(|x|^{2n +1}+|y|^{2n+1}), \text{ \ as $x,y \rightarrow  0$}.$$
As in the last subsection, we still let the $\Torus$-directional $2$-jets of $\Phi$ along $\ZentralFaser$ be
\[\aoneone x^2+\atwotwo y^2,\]
and denote $$\widetilde L=L-\aoneone x^2-\atwotwo y^2.$$

\begin{notation}\label{SpaceofHamogenousPolynomials}
In the following, we denote the space of  homogenous polynomials in $x,y$ with real coefficients and of degree $2n$ by $\HappytoN$. In $\HappytoN$, the polynomials which are even in both variables form a linear subspace. We denote this subspace by $\HappytoNE$. 

\end{notation}
Because of the assumption that $\Phi$ is even in both $x$ and $y$, we have that $P$ is also even in both $x$ and $y$. In the following, we view $P(t,\ast)$ as a curve in $\HappytoNE$ and try to derive an ODE for $P$.
 To do this, we plug (\ref{RevisedTaylorDec302018}) into (\ref{Gleichung:GeodesicEquationonToruswithFlatBackGround}) and concentrate on terms of degree $2n$. We assumed that $\Phi\in C^{2n+4}$, so $R'$ and $R''$ are both controlled by $O(|x|^{2n+1}+|y|^{2n+1})$ as $x,y\rightarrow 0.$ Then using the fact that 
 the lowest terms in $\Phi $, $\Phi'$ and $ \Phi''$, if not zero, are of degree $2$, we get that, in (\ref{Gleichung:GeodesicEquationonToruswithFlatBackGround}), $R $ does not contribute to any terms of degree $2n$. And using that the lowest terms in $\widetilde L $,  $\widetilde L' $ and $\widetilde L''$, if not zero, are of degree greater than $2$, we find that all terms in 
\[\widetilde L'' \cdot \Lappen P, \ \ P''\cdot \Lappen \widetilde L,\ \  \Naben P'\cdot \Naben \widetilde L' \]
are of degree greater than $2n$, so they also do not contribute to terms of degree $2n$ in (\ref{Gleichung:GeodesicEquationonToruswithFlatBackGround}).  After eliminating these terms, we find that the $2n$th degree part of (\ref{Gleichung:GeodesicEquationonToruswithFlatBackGround}) is 
\begin{equation}(1+2\aoneone+2\atwotwo)P''+(a''x^2+b''y^2)\Lappen P-4(\aoneone' x, \atwotwo' y)\cdot\Naben P'- K_1(L)=0.
																			\label{gleichung:2nthordertermDec302018}
\end{equation}
Above, $K_1(L)$ is the $2n$th degree part of 
\[-(\Lappen L)L''+|\Naben L'|^2,\]
so it is determined by $L$. We divide (\ref{gleichung:2nthordertermDec302018}) by $1+2\aoneone+2\atwotwo$, and use relation (\ref{Gleichung:Derivativeofa11a22overVolumeform}) (\ref{Gleichung:secondDerivativeofa11a22overVolumeform}), then  (\ref{gleichung:2nthordertermDec302018}) becomes 
\begin{equation}\label{gleichung:fortgeschrittenODEP}
P''+\left(4\epsilon^2 A^2x^2+\frac{4\epsilon^2}{A^2}y^2\right)\Lappen P-4\epsilon\left(Ax, -\frac{1}{A}y\right)\cdot\Naben P'=\frac{K_1(L)}{1+2\aoneone+2\atwotwo}.
\end{equation}
We denote
\begin{equation}\label{definitionofMEAMSA}
\left( A^2x^2+\frac{1}{A^2}y^2\right)\Lappen=\MEA,\ \ \ \ \ \left(Ax, -\frac{1}{A}y\right)\cdot\Naben=\MSA,
\end{equation}
and consider them as elements of $\Endo(\HappytoNE)$.

To simplify (\ref{gleichung:fortgeschrittenODEP}), a standard idea is  to find a family of automorphisms of $\HappytoNE$,  $\MU(t)$: 
$[0,1] \rightarrow\Auto$($\HappytoNE$), and let 
\[P(t)=\MU(t)Q(t),\]
so that the homogenous part of the equation for $Q$ becomes 
\[Q''(t)=T(t)Q(t),\]
for some $T: [0,1]\rightarrow \Endo(\HappytoNE)$. In the following, we show that this can be done  and it turns out that $T(t)$ is independent of $t$.

Plug $P(t)=\MU(t)Q(t)$ into (\ref{gleichung:fortgeschrittenODEP}), and let $\MU^{-1}$ act on  (\ref{gleichung:fortgeschrittenODEP}) from left, to get
\begin{equation}\label{SecondorderODEofQ}
Q''+2 \MU^{-1}(\MU'-2\epsilon \MSA\MU)Q'+\MU^{-1}(\MU''+4\epsilon^2 \MEA \MU-4\epsilon \MSA \MU')Q=\MU^{-1}\frac{K_1}{1+2\aoneone+2\atwotwo}.
\end{equation}
To make
\begin{equation}\MU'-2\epsilon \MSA \MU=0,	\label{ODEofMUDec302018}
\end{equation}
we can choose, for some constant $C_0$, 
\begin{equation}\label{SymbolicalExpressionofU}
\MU(t)=\exp\left(2\epsilon\int^t_0 \MSA+C_0\right)=\exp\left(2\epsilon\int^t_0 A x\partial_x-2\epsilon\int_0^t\frac{1}{A}y\partial y+{C_0}\right).
\end{equation}
Equation (\ref{ODE_of_A_hyperbolic}) gives that
\[2\epsilon \int^t_0 A=\frac{\log(A^2+1)}{2}+C_1, \ \ \ \ \  -2\epsilon\int_0^t \frac{1}{A}={1\over 2}\log\left(1+\frac{1}{A^2}\right)+C_2\]
Plugging these into (\ref{SymbolicalExpressionofU}), we find that we can let
\[\MU(t)=\exp\left(\frac{\log(A^2(t)+1)}{2}x\partial_x\right) \exp\left(\frac{\log(1+\frac{1}{A^2(t)})}{2}y\partial_y\right),\] since $x\frac{\partial }{\partial x}$ commutes with $y\frac{\partial }{\partial y}$ and (\ref{ODEofMUDec302018}) is linear, so a constant multiple of $\MU$ also satisfies (\ref{ODEofMUDec302018}).
All monomials $x^j y^k$ are eigenvectors of $U$,  with eigenvalues
\[(A^2+1)^{\frac{j}{2}}(A^{-2}+1)^{\frac{k}{2}},\] since $x^j y^k$ are eigenvectors of $x\frac{\partial }{\partial x}$ and $y\frac{\partial }{\partial y}$, with eigenvalues $j$ and $k$.

Differentiating (\ref{ODEofMUDec302018}) with respect to $t$, gives 
\begin{equation}\label{secondorderderivativeofMU}
\MU''=4\epsilon^2\MSA \MSA \MU+2\epsilon \MSA' \MU.
\end{equation}
Differentiating the definition of $\MSA$ in (\ref{definitionofMEAMSA}) and using (\ref{ODE_of_A_hyperbolic}), gives
\begin{equation}\label{DerivativeofMSA}
\MSA'=2\epsilon\left[(A^2+1)x\partial_x+\left(1+\frac{1}{A^2}\right)y\partial_y\right].
\end{equation}

Then plugging (\ref{secondorderderivativeofMU})(\ref{DerivativeofMSA}) into (\ref{SecondorderODEofQ}) and using direct algebraic computation,  (\ref{SecondorderODEofQ}) becomes
\begin{equation}\label{simplifiedSecondorderODEofQ}
Q''+4\epsilon^2\MU^{-1}((A^2+1)x\partial_x+(1+{1\over A^2})y\partial_y+\MEA-\MSA\MSA)\MU Q=\MU^{-1}\frac{K_1}{1+2\aoneone+2\atwotwo}.
\end{equation}
To simplify (\ref{simplifiedSecondorderODEofQ}), we use (\ref{definitionofMEAMSA}) and the following explicit expression of $\MSA\MSA$, 
\[\MSA \MSA=A^2x^2\partial^2_x+\frac{y^2}{A^2}\partial_y^2-2xy\partial_x\partial_y+A^2x\partial_x+\frac{y}{A^2}\partial_y.\]
This gives
\begin{align}
\MU^{-1}\frac{K_1}{1+2\aoneone+2\atwotwo}&=Q''+4\epsilon^2\MU^{-1}(A^2x^2\partial_y^2+\frac{1}{A^2}y^2\partial_x^2+2xy\partial_x\partial_y+x\partial_x+y\partial_y)\MU Q\\
&=Q''+4\epsilon^2\MU^{-1}(Ax\partial_y+\frac{1}{A}y\partial_x)^2\MU Q\\
&=Q''+4\epsilon^2\left[\MU^{-1}(Ax\partial_y+\frac{1}{A}y\partial_x)\MU\right]^2 Q.     \label{FinalQ}
\end{align}
The most crucial computation of this section is that, in (\ref{FinalQ}),
\begin{equation}\label{firstappearBoostOperator}
\MU^{-1} (Ax\partial_y+{1\over A}y\partial_x) \MU=x\partial_y+y\partial_x,
\end{equation}
which says that $\MU^{-1} (Ax\partial_y+{1\over A}y\partial_x) \MU$ is independent of $t$. (\ref{firstappearBoostOperator}) can be verified by letting operators on both sides of it act on  $x^\mu y^\nu\in\HappytoN$, for $\mu,\nu\in \EZ^{\geq 0}$. Note that $x\partial_y+y\partial_x$ and $Ax\partial_y+{1\over A}y\partial_x$ are not elements of $\Endo(\HappytoNE)$, instead they are elements of $\Endo(\HappytoN)$, however, their squares are elements of $\Endo(\HappytoNE)$.

We denote
\[x\partial_y+y\partial_x=\Boost,\]
then the  equation of $Q$ becomes
\begin{equation}\label{NewNewNewEquationofQ}
Q''+4\epsilon^2\Boost^2Q=\MU^{-1}\frac{K_1(L)}{1+2\aoneone+2\atwotwo},
\end{equation}and the right hand side is determined by $L$.

It's easy to check that
\begin{equation}q_k=(x+y)^{n+k}(x-y)^{n-k}+(x+y)^{n-k}(x-y)^{n+k}, \ \ \ \text{for $k=0,\ ...,\  n$},\label{RefereeEigenBasis}
\end{equation}
are eigenvectors of $\Boost^2$ in $\HappytoNE$, with eigenvalues $(2k)^2$, and so form an eigenbasis.

Then using basic spectral theory, we know if $4\epsilon^2\cdot (2n)^2<\pi^2$, $Q$ is determined by its boundary value and $L$. This means that along $[0,1]\times \{0\}$, the $\Torus$-directional $2n$-jets of $\Phi$ is determined by the $\Torus$-directional $(2n-2)$-jets of $\Phi$ and the $2n$-jets of $\varphi_0$ and $\varphi_1$ at $0$. By induction and Theorem \ref{Satz:SecondJetsLemmaA}, we know the $\Torus$-directional $2n-$jets of $\Phi$ is determined by the $2n-$jets of $\varphi_0$ and $\varphi_1$ at $0$.

If $4\epsilon^2\cdot (2n)^2=\pi^2$, we can still use the arguments above to show the $\Torus$-directional $(2n-2)-$jets of $\Phi$ is determined by the $(2n-2)-$jets of $\varphi_0$ and $\varphi_1$ at $0$. So, in equation (\ref{NewNewNewEquationofQ}), the right hand side is determined. However, to make (\ref{NewNewNewEquationofQ}) solvable, $Q(0)$ and $Q(1)$ have to satisfy some compatibility conditions, which we explain in the following. 

In (\ref{NewNewNewEquationofQ}) we let 
\[Q=\sum_{j=0}^n f_j q_j, \ \ \ \ \ \ \MU^{-1}\frac{K_1(L)}{1+2\aoneone+2\atwotwo}=\sum_{j=0}^n k_j q_j,\]
where $f_i$'s and $k_i$'s are functions of $t$ defined on $[0,1]$, and $q_i$'s are given by (\ref{RefereeEigenBasis}).
Then the equation for $Q$ becomes
\[0=(f_0''-k_0)q_0+...+(f_n''+4\epsilon^2\cdot(2n)^2f_n-k_n)q_n.\]
So $f_n$ must satisfy
\begin{equation}
f_n''+\pi^2f_n-k_n=0.\label{equationoffnDec302018}
\end{equation}
Multiplying (\ref{equationoffnDec302018}) by $\sin(\pi t)$ and integrating on $[0,1]$, gives
\begin{equation}
				\label{CompatibleConditionoffn}
f_n(0)+f_n(1)=-f_n\cdot \cos(\pi t)\big |_0^1=\int_0^1\frac{k_n(t)\sin(\pi t)}{\pi} dt.
\end{equation}

We want to transform (\ref{CompatibleConditionoffn}) into a condition of $\varphi_0$ and $\varphi_1$. Since $q_j$, $j=0,\ ...\ ,n$,  are linearly independent, we can find a $2n$th order partial differential operator 
$$\MtD=\sum_{j=0}^n\nu_j\partial_x^{2j}\partial_y^{2n-2j},$$ 
so that
\[\MtD(q_j)=0, \text{ for } j=0,\ ...\ ,n-1,\ \  \text{and\ }\MtD(q_n)=1.\]
Now we have $\MtD(Q)=f_n$,  then we try to modify $\MtD$ to get an operator $\MD$, so that, $\MD(P)=f_n$. To do this
  we want \[\MD(\MU g)=\MtD(g),   \ \ \ \text{for any }g\in \HappytoNE,\] so we can let
$$\MD=\sum_{j=0}^n\nu_j\frac{A^{2n-2j}}{(1+A^2)^n}\partial_x^{2j}\partial_y^{2n-2j},$$
and get 
\[\MD(\Phi)|_{\ZentralFaser}=\MD(P)=f_n.\]

With the operator $\MD$, we can write (\ref{CompatibleConditionoffn}) as
\[\MD\varphi_0+\MD\varphi_1=K,\]
where $K$ is determined by the $2n-2$ jets of $\varphi_0$ and $\varphi_1$ at $0$.

Now, we reach the following conclusion of this section
\begin{Satz}\label{Satz:Paaren}
For $n\in\EZ$, $n\geq 2$, given $\Phi\in C^{2n+\Add}([0,1]\times\Torus)$, a non-degenerate geodesic which is even in both $x$ and $y$
and satisfies
\[\Phi(0,\ast)=\varphi_0,\ \ \ \ \Phi(1,\ast)=\varphi_1,\]
 we have 
\begin{equation}\frac{\Phi_{xx}'\Phi_{yy}'}{4(1+\Phi_{xx}+\Phi_{yy})^2}\label{ANegativeConstant}
\end{equation}
is a constant along $\ZentralFaser$. It is a negative constant if and only if 
\[(\partial_{xx}\varphi_1-\partial_{xx}\varphi_0)(0)\cdot(\partial_{yy}\varphi_1-\partial_{yy}\varphi_0)(0)<0.\] 
In this case, denoting the negative constant by $-\epsilon^2$, for some positive constant $\epsilon$,
we have that $\epsilon$ is determined by the $2$-jets of $\varphi_0$ and $\varphi_1$. And
 
(1) if 
\begin{equation}n^2\cdot 16\epsilon^2< \pi^2,\label{SoftSituationHighDerivative}
\end{equation}
 along $\ZentralFaser$ all  the derivatives of $\Phi$ of order less than or equal to $2n$ are determined by the $2n$-jets of $\varphi_0,\varphi_1$ at $0\in\Torus$;

(2) if
\begin{equation}n^2\cdot 16\epsilon^2= \pi^2,\label{RigidSituationHighDerivative}
\end{equation}
then, along $\ZentralFaser$ all the derivatives of $\Phi$ of order less than $2n$ are determined by the $(2n-2)$-jets of $\varphi_0,\varphi_1$ at $0\in\Torus$; and there are two non-zero vectors, $\vu$ and $ \vv$,
so that $\varphi_0$ and $\varphi_1$ have to satisfy
\begin{equation}
\sum_{i=0}^n v_i\cdot{D_x^{2n-2i}}{D_y^{2i}\varphi_0}(0)+\sum_{i=0}^n u_i\cdot{D_x^{2n-2i}}{D_y^{2i}\varphi_1}(0)
=K,																						\label{Confundido}
\end{equation}
where $\vu,\ \vv$ and $K$ are all determined by the $(2n-2)$-jets of $\varphi_0,\varphi_1$ at $0\in\Torus$.
\end{Satz}

\section{Construction}\label{KurzUnGlatt}
In this section we construct the example claimed in Example \ref{Beispiel:Kurz} of Section \ref{IntroductionSection}. The main idea is to choose $\varphi_1$ and $\varphi_0$, so that our geodesic falls into the case (2) of Theorem \ref{Satz:Paaren} and (\ref{Confundido}) is violated. 

Given $n\in \EZ$, $n>2$, suppose we can connect $0$ and 
\[h_n=\frac{1}{2}\sin \left(\frac{\pi}{2n}\right)\cdot(\sin^2x-\sin^2y),\]
by a $C^{2n+\Add}$ geodesic $\Phi$,   then according to (\ref{LemmaA_Formula_for_epsilon}) of Theorem \ref{Satz:SecondJetsLemmaA} and the analysis of Section \ref{ODESystemofSecondOrderDerivativesofGeodesic}, we have 
that along $\ZentralFaser$
\[\frac{\Phi_{xx}'\Phi_{yy}'}{4(1+\Phi_{xx}+\Phi_{yy})^2}=-\frac{\pi^2}{16n^2}.\]

Then using Theorem \ref{Satz:Paaren}, we can find $\vv\in \ER^{n+1}$, with $v_\kappa\neq 0$, for some $\kappa\in \EZ$,
s.t. 
\begin{equation}
\sum_{i=0}^n v_i\cdot{D_x^{2n-2i}}{D_y^{2i}h_n}(0)
=K.	\label{KurzPaaren}
\end{equation}
Above, $\vv$ and $K$ only depend on the $(2n-2)$-jet of $h_n$ at $0\in\Torus$.

Then for any $\chi\in \ER, \chi\neq 0$,  if 
\begin{equation}
\widetilde h_n=\frac{1}{2}\sin \left(\frac{\pi}{2n}\right)\cdot(\sin^2x-\sin^2y)+\chi \sin^{2n-2\kappa}x\sin^{2\kappa}y,
\end{equation}
can also be connected with
$0$ by a $C^{2n+\Add}$ non-degenerate geodesic, using Theorem \ref{Satz:Paaren} again, we have
\begin{equation}
\sum_{i=0}^n u_i\cdot D_x^{2n-2i}D_y^{2i}\widetilde h_n(0)
=J,	\label{KurzPaaren_Added}
\end{equation} 
with $\bf u$ and $J$ only depending on the $(2n-2)$-jets of $\widetilde h_n$ at $0\in\Torus$.
Note that $\widetilde h_n$ and $h_n$ have the same $(2n-2)$-jets at $0\in\Torus$, so the $\bf v$(or $K$)  in (\ref{KurzPaaren}) is equal to the $\bf u$(or $J$)  in (\ref{KurzPaaren_Added}). Taking difference of (\ref{KurzPaaren}) and (\ref{KurzPaaren_Added}), we  get
\[v_\kappa\cdot(2n-2\kappa)!(2\kappa)!=0,\]which is a contradiction. 

Now let the $\chi$ be $e^{-n}$, we find for any $n\in\EZ$, $n\geq 3$, either
\[h_n=\frac{1}{2}\sin \left(\frac{\pi}{2n}\right)\cdot(\sin^2x-\sin^2y),\]or 
\[\widetilde{h}_n=\frac{1}{2}\sin \left(\frac{\pi}{2n}\right)\cdot(\sin^2x-\sin^2y)+e^{-n} \sin^{2n-2\kappa}x\sin^{2\kappa}y,\]
cannot be connected with
$0$ by a non-degenerate $C^{2n+\Add}$ geodesic. 

It is easy to see, for any fixed $B>0$, $|h_n|_B+|\widetilde h_n|_B\rightarrow 0$ as $n\rightarrow \infty.$ So we can pick a subsequence from  $\{h_n\}_{n=3}^\infty\cup\{\widetilde h_n\}_{n=3}^\infty$ satisfying the requirement of Example \ref{Beispiel:Kurz}.

\appendix
\section{The Boundary Value Problem for an ODE System}\label{Solvability_of_2nd_Jets}
In this appendix, we discuss the following problem

\begin{problem}\label{Problem:2JetsODE}
Given real numbers: $\aoneone_0 $, $\atwotwo_0 $, $\aoneone_1 $ and $\atwotwo_1 $, satisfying 
\begin{align}\renewcommand\arraystretch{1.3}
\aoneone_1+\atwotwo_1+{1\over 2}>0\label{11positive};\\
\aoneone_0+\atwotwo_0+{1\over 2}>0\label{00positive},
\end{align}
 find $\aoneone(t) $, $\atwotwo(t) \in C^2([0,1];\ER)$, \  satisfying the following equations\ \footnote{Denoting $\frac{d}{dt}$ by prime ``\ $'$\ ".}
\begin{align}\renewcommand\arraystretch{1.3}
\aoneone ''=\frac{4(\aoneone ')^2}{1+2\aoneone +2\atwotwo }		;		\label{appendixEvolutionofa11}\\
\atwotwo ''=\frac{4(\atwotwo ')^2}{1+2\aoneone +2\atwotwo }	,	\label{appendixEvolutionofa22}
\end{align}
the boundary conditions
\[\aoneone(0)=\aoneone_0,\ \ \aoneone(1)=\aoneone_1,\ \ \atwotwo(0)=\atwotwo_0,\ \ \atwotwo(1)=\atwotwo_1,\]
and
\begin{equation}
\aoneone +\atwotwo +{1\over 2}>0,  \ \ \ \ \text{ on\ \ \ \  }[0,1],
\end{equation}

\end{problem}
We will prove the following theorem, regarding the solvability of Problem \ref{Problem:2JetsODE} and some properties of the solution.

\begin{Satz}{
}\label{Satz:SecondJetsLemmaA}
The necessary and sufficient condition for Problem  \ref{Problem:2JetsODE} to have a solution is that the boundary data satisfy
\begin{align}\renewcommand\arraystretch{1.3}
\aoneone_0+\atwotwo_1+{1\over 2}>0\label{01positive};\\
\aoneone_1+\atwotwo_0+{1\over 2}>0\label{10positive}.
\end{align}
The solution is then uniquely determined by the boundary data.  For a solution $(\aoneone, \atwotwo)$,
\begin{equation}
\frac{\aoneone '}{1+2 \aoneone +2 \atwotwo }\cdot \frac{\atwotwo '}{1+2 \aoneone +2 \atwotwo }\label{Dec28ProductofEigenValue}
\end{equation}
is a constant on $[0,1]$, it's a negative constant if  and only if
\begin{equation}(\aoneone_0-\aoneone_1)\cdot(\atwotwo_0-\atwotwo_1)<0,       \label{DifferentSignConditionLemmaA}
\end{equation}
And in this case, we denote (\ref{Dec28ProductofEigenValue}) by $-\epsilon^2$, for some positive number $\epsilon$. 
We have that 
\[\epsilon<{\pi\over 4},\] and $\epsilon$ is determined by
\begin{equation}\label{LemmaA_Formula_for_epsilon}
\cos(4\epsilon)=\frac{(1+2\aoneone(0)+2\atwotwo(0))^2+(1+2\aoneone(1)+2\atwotwo(1))^2-(2\aoneone(1)-2\atwotwo(1)-2\aoneone(0)+2\atwotwo(0))^2}
{2(1+2\aoneone(0)+2\atwotwo(0))(1+2\aoneone(1)+2\atwotwo(1))}.
\end{equation}
\end{Satz}

We will transform Problem \ref{Problem:2JetsODE} to the problem of finding geodesics in a constant curvature Lorentz manifold. 

Let us denote 
\begin{equation}
1+2\aoneone+2\atwotwo=Z;
\end{equation}
\begin{equation}
2\aoneone-2\atwotwo=X.
\end{equation}

Then (\ref{appendixEvolutionofa11})(\ref{appendixEvolutionofa22}) become 
\begin{equation}
Z''=\frac{Z'^2+X'^2}{Z};     \label{GeodesicEquationforZeit}
\end{equation}
\begin{equation}
X''=\frac{2Z'X'}{Z}.           \label{GeodesicEquationforSpace}
\end{equation}

We find that the ODE system above is actually the Euler-Lagrange equation of the functional
\[\int_0^1\frac{X'^2-Z'^2}{Z^2}dt.\]
This implies that the ODE system  (\ref{GeodesicEquationforZeit}) (\ref{GeodesicEquationforSpace}) is actually the geodesic equation for the upper half-plane
\[\LP\triangleq \{(X,Z)\big|X\in \ER,\ Z\in \ER^{+}\},\]
 with the metric
\begin{equation}ds^2=\frac{dX^2-dZ^2}{Z^2}.    \label{LorentzPoincareMetric}
\end{equation}

In \cite{KatsumiNomizu}, the metric (\ref{LorentzPoincareMetric}) is referred to as the Lorentz-Poincar\'e metric. If we let $u=\log(Z)$, then
$(X,u)$ are called the planar(or flat slicing) coordinates of de Sitter space \cite{LesHouches} \cite{Hartman}.

Let $(\aoneone(t), \atwotwo(t))$ be a solution to Problem \ref{Problem:2JetsODE}, then 
\[\gamma(t)\triangleq (X(t), Z(t))=(2\aoneone(t)-2\atwotwo(t), 1+2\aoneone(t)+2\atwotwo(t)),\] is a geodesic\footnote{In this appendix, a parametrized curve is called a geodesic if it satisfies the geodesic equation.} in $\LP$. Since 
\begin{equation}
-\frac{\aoneone'\atwotwo'}{(1+2\aoneone+2\atwotwo)^2}=-\frac{(X'+Z')(Z'-X')}{16Z^2}=\frac{X'^2-Z'^2}{16Z^2}=\frac{ds^2(\gamma',\gamma')}{16},
\end{equation}
the fact that $\frac{\aoneone'\atwotwo'}{(1+2\aoneone+2\atwotwo)^2}$ is a constant
\footnote{ This is also proved by (\ref{InvarianceofBifurcatingStationaryLocus}).} 
is equivalent to the fact  that along the geodesic $\gamma$, the kinetic energy $ds^2(\gamma',\gamma')$ is a constant.

The sign of $\frac{\aoneone'\atwotwo'}{(1+2\aoneone+2\atwotwo)^2}$ is related to the casual character of $\gamma$, providing $\gamma'\neq 0$, as shown in the following:
\begin{itemize}
\item$\frac{\aoneone'\atwotwo'}{(1+2\aoneone+2\atwotwo)^2}>0$ $\Leftrightarrow$ $X'^2-Z'^2<0$ $\Leftrightarrow$ $\gamma$ is a time-like geodesic;
\item$\frac{\aoneone'\atwotwo'}{(1+2\aoneone+2\atwotwo)^2}<0$ $\Leftrightarrow$ $X'^2-Z'^2>0$ $\Leftrightarrow$ $\gamma$ is a space-like geodesic;
\item$\frac{\aoneone'\atwotwo'}{(1+2\aoneone+2\atwotwo)^2}=0$ $\Leftrightarrow$ $X'^2-Z'^2=0$ $\Leftrightarrow$ $\gamma$ is a light-like geodesic.
\end{itemize}

Suppose that two points $(X_0, Z_0)$ and $(X_1, Z_1)$ are connected by a space-like geodesic $\gamma(t)=(X(t), Z(t))$, for $t\in[0,1]$, then we have
\[\epsilon\triangleq\sqrt{-\frac{\aoneone'\atwotwo'}{(1+2\aoneone+2\atwotwo)^2}}=\frac{1}{4}\sqrt{ds^2(\gamma',\gamma')}=\int_0^1\frac{1}{4}\sqrt{ds^2(\gamma',\gamma')}.\]
So $4\epsilon$ is simply the length of the geodesic.  

As shown in \cite{KatsumiNomizu}, all geodesics in $\LP$ are either straight lines
\[\{X\equiv X(0)\},\]
or arcs of the following hyperbolas
\[\{(X, Z)\big| Z^2-(X-\lambda)^2=C\},\]
for constants $b\in \ER, C\in \ER$.\footnote{$\{X\equiv X(0)\}$ is a time-like geodesic. The hyperbola is space-like geodesic if $C>0$; it is time-like geodesic if $C<0$; it is light-like if $C=0$. When $C=0$, the hyperbola degenerates and it is the union of two straight lines.} Given two points $(X_0, Z_0)$ and $(X_1, Z_1)$, if $X_0\neq X_1$, we can uniquely determine the hyperbola passing through them, if $X_0= X_1$ the geodesic connecting them must be a straight line, so there can be at most one   geodesic connecting any two points. This implies that $4\epsilon$ is also the geodesic distance between $(X_0, Z_0)$ and $(X_1, Z_1)$.
When computing the hyperbola passing through two points, if we get a hyperbola that is not connected in $\LP$, corresponding to $C\leq 0$, then the two points may be located on two different  branches of the hyperbola, and in this case we don't say that these two points can be connected by a geodesic.

With the observations above, we can equivalently translate Theorem \ref{Satz:SecondJetsLemmaA} to the following Theorem \ref{Satz:Geodesics_deSitter_LemmaAA}

\begin{Satz}\label{Satz:Geodesics_deSitter_LemmaAA}
Two points $(X_0, Z_0)$ and $(X_1, Z_1)$ in $\LP$ can be connected by a geodesic if and only if
\begin{equation}
Z_1+Z_0>|X_1-X_0|.			\label{ConnectableConditionLemmaAA}		
\end{equation}
And when the condition above is satisfied, there is only one geodesic connecting these two points.
If, in addition, 
\begin{equation}
|Z_1-Z_0|<|X_1-X_0|,			\label{SpacelikeConnectableConditionLemmaAA}		
\end{equation}
then, $(X_0, Z_0)$ and $(X_1, Z_1)$ can be connected by a space-like geodesic, their geodesic distance 
\[D\triangleq D((X_0, Z_0),(X_1, Z_1))\]satisfies
\begin{equation}
D<\pi,				\label{LengthLimitofSpacelikeGeodesic2019}
\end{equation}
and is determined by 
\begin{equation}
\cos(D)=\frac{Z_0^2+Z_1^2-(X_0-X_1)^2}{2Z_0Z_1}.		\label{ArcLengthFormulaLemmaAA}		
\end{equation}
\end{Satz}

The transition from Theorem \ref{Satz:SecondJetsLemmaA} to Theorem \ref{Satz:Geodesics_deSitter_LemmaAA} is straightforward. What we need to do is simply 
replacing $1+2\aoneone+2\atwotwo$ by $Z$, 
replacing $2\aoneone-2\atwotwo$ by $X$,
replacing ``solution to Problem \ref{Problem:2JetsODE}" by ``geodesic in $\LP$'',
replacing $4\epsilon$ by $D$,
and using the fact that if $a'b'<0$, the corresponding geodesic in $\LP$ is space-like. Then we find conditions  (\ref{01positive}) (\ref{10positive}) correspond to (\ref{ConnectableConditionLemmaAA}), the condition (\ref{DifferentSignConditionLemmaA}) corresponds to (\ref{SpacelikeConnectableConditionLemmaAA}), 
and the relation (\ref{LemmaA_Formula_for_epsilon}) exactly becomes (\ref{ArcLengthFormulaLemmaAA}). So, if Theorem \ref{Satz:Geodesics_deSitter_LemmaAA} is proved, Theorem \ref{Satz:SecondJetsLemmaA} follows immediately.
In the following, we prove Theorem \ref{Satz:Geodesics_deSitter_LemmaAA}.

\begin{proof}

Using high school computation, we can show that all geodesics passing through $(0,1)$ cover the area 
\[\{(X,Z)\big|Z+1>|X|\},\]
and nothing more.
As illustrated by Figure \ref{fig:NewPicture_Hyperbolas}, the gray area $\{Z+1>|X|\}\cap\{|Z-1|<|X|\}$ is covered by the following family of space-like geodesics $$\{(X,Z)\big|Z^2-(X-\lambda)^2=1-\lambda^2\},\ \ \ \text{for } -1<\lambda<1.$$ The dotted area $\{|Z-1|>|X|\}$ is covered by the following family of time-like geodesics:
\begin{itemize}
\item the left-branch of $\{(X,Z)\big|Z^2-(X-\lambda)^2=1-\lambda^2\},\ \ \ \text{for } 1<\lambda$;
\item the right-branch of $\{(X,Z)\big|Z^2-(X-\lambda)^2=1-\lambda^2\},\ \ \ \text{for } \lambda<-1$;
\item $\{X\equiv 0\}.$
\end{itemize}
Further, $\{|Z-1|=|X|\}$ is the union of two light-like geodesics.
\begin{figure}
\centering
\includegraphics[height=4.2cm]{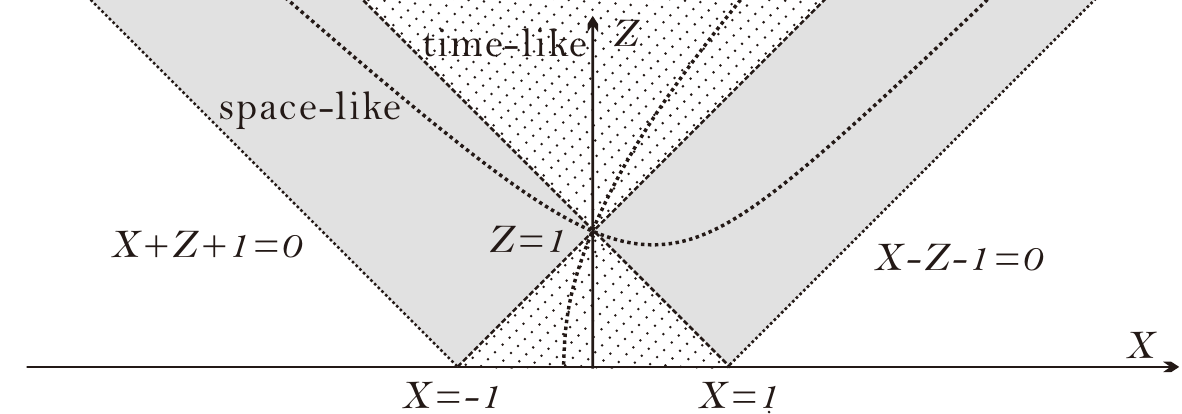}
\caption{Area covered by hyperbolas passing through $(0,1)$}
\label{fig:NewPicture_Hyperbolas}
\end{figure}

Now given two points $(X_0, Z_0)$ and $(X_1, Z_1)$, we consider the following isometric transformation on $\LP$,
\[\ReducetoBase(X,Z)=\left(\frac{X-X_0}{Z_0}, \frac{Z}{Z_0}\right).\]
Because $\ReducetoBase$ is an isometry, $(X_0, Z_0)$ and $(X_1, Z_1)$ can be connected by a geodesic if and only if $\ReducetoBase(X_0, Z_0)$ and $\ReducetoBase(X_1, Z_1)$ can be connected by a geodesic. This is equivalent to
\begin{equation}\ReducetoBase(X_1, Z_1)\in\{(X,Z)\big|Z+1>|X|\}, \label{StaysinShadowedandDottedAreaLemmaAA}
\end{equation}
since $\ReducetoBase(X_0, Z_0)=(0,1).$
Obviously, (\ref{StaysinShadowedandDottedAreaLemmaAA}) is equivalent to
\[\frac{Z_1}{Z_0}+1>\left|\frac{X_1-X_0}{Z_0}\right|\]
and (\ref{ConnectableConditionLemmaAA}), since $Z_0>0$.

Similarly, the geodesic connecting $(X_0, Z_0)$ and $(X_1, Z_1)$ is a space-like geodesic if, in addition to (\ref{StaysinShadowedandDottedAreaLemmaAA}), 
\begin{equation}\ReducetoBase(X_1, Z_1)\in\{(X,Z)\big||Z-1|<|X|\}.\label{StayinShadowedArea}
\end{equation}
And (\ref{StayinShadowedArea}) is equivalent to 
\[\left|\frac{Z_1}{Z_0}-1\right|<\left|\frac{X_1-X_0}{Z_0}\right|\]
and (\ref{SpacelikeConnectableConditionLemmaAA}).

Using the formula below (9') of \cite{KatsumiNomizu}, we know the length of any space-like geodesic in $\LP$ is smaller than $\pi$, which proves (\ref{LengthLimitofSpacelikeGeodesic2019}). This also enables us to use formula (33) of \cite{LesHouches}\footnote{Note that $(X,\log(Z))$ are the planar coordinates for de Sitter space.} and get
\[\cos(D)=\cosh(\log Z_1-\log Z_0)-\frac{1}{2}\exp\left({-\log Z_1-\log Z_0}\right)(X_1-X_0)^2=\frac{Z_1^2+Z_0^2-(X_1-X_0)^2}{2Z_1Z_2},\]
which is exactly formula (\ref{ArcLengthFormulaLemmaAA}). Note that conditions (\ref{ConnectableConditionLemmaAA}) and (\ref{SpacelikeConnectableConditionLemmaAA}) guarantee that the right hand side of (\ref{ArcLengthFormulaLemmaAA}) is smaller than $1$ and greater than $-1$.
\end{proof}

\section*{Acknowledgement}
The author is partially supported by National Natural Science Foundation of China (grant no. 11571330, 11271343). He would like to thank  Prof. Xiuxiong Chen, Prof. Mikhail Feldman,  Jiyuan Han, Jingrui Cheng  and Long Li for very helpful discussion. He also wants to thank the anonymous referee for his/her many insightful and thoughtful suggestions on improving the paper.


Jingchen Hu\\
Institute of Mathematical Sciences, ShanghaiTech University. 393 Middle Huaxia Road, Shanghai, 201210, China.\\
Email:
JINGCHENHOO@GMAIL.COM

\end{document}